\documentclass[a4paper,10pt]{amsart}

\usepackage{amsmath,amsfonts,amsthm, amssymb,xypic, epsfig}
\usepackage[all]{xy}

\def\C{\mathbb{C}}
\def\Z{\mathbb{Z}}
\def\N{\mathbb{N}}
\def\qed{$\hfill \checkmark$}
\def\U{\mathbf{U}}
\def\H{\mathbf{H}}

\def\g{\mathfrak{g}}
\def\n{\mathfrak{n}}

\def\O{{\mathcal{O}}}

\def\m{\mathbf{m}}
\def\d{\mathbf{d}}
\def\A{\mathbb{A}}

\def\d{\mathbf{d}}
\def\i{\mathbf{i}}
\def\j{\mathbf{j}}
\def\k{\mathbf{k}}
\def\barv{\overline{\phantom{x}}}
\def\fq{\overline{\mathbb{F}_q}}

\newtheorem{theo}{Theorem}
\newtheorem{prop}{Proposition}
\newtheorem{lem}{Lemma}
\newtheorem{cor}{Corollary}

\numberwithin{equation}{section}
\numberwithin{lem}{section}
\numberwithin{prop}{section}

\numberwithin{cor}{section}
\numberwithin{theo}{section}
\numberwithin{claim}{section}

\begin{document}

\title[Canonical bases for quantum generalized Kac-Moody algebras]
{Canonical bases for Quantum \\
Generalized Kac-Moody algebras}

\author{Seok-Jin Kang$^*$ and Olivier Schiffmann$\dag$}
\date{}

\address{$^*$ School of Mathematics\\
          Korea Institute for Advanced Study\\
          207-43 Cheongryangri-Dong\\
          Dongdaemun-Gu\\
          Seoul 130-722, Korea}
\email{sjkang@kias.re.kr}
\address{$^\dag$ DMA, \'Ecole Normale Sup\'erieure, 45 rue d'Ulm,
          75230 Paris Cedex 05-FRANCE}
\email{schiffma@dma.ens.fr}

\thanks{$^*$This research was supported by KOSEF Grant
\# R01-2003-000-10012-0 and KRF Grant \# 2003-070-C00001}

\begin{abstract}

We construct canonical bases for quantum generalized Kac-Moody algebras
using semisimple perverse sheaves.

\end{abstract}

\maketitle

\section*{Introduction}

\vspace{.1in}

The {\it quantum groups} arise as certain deformations of
the universal enveloping algebras of Kac-Moody algebras.
One of the main achievements in the theory of quantum groups
is the construction of {\it canonical bases}
(or {\it global bases}) for integrable highest weight modules
over a quantum group $\U_v(\g)$.
The canonical bases enjoy many remarkable properties
such as positivity with respect to the $\U_v(\g)$-action,
and they all come from a single canonical basis for
the negative part $\U_v^-(\g)$ of the quantu group
by acting on the highest weight vector.
These bases were constructed independently
by Kashiwara (in a purely algebraic way) and Lusztig
(in a geometric way via intersection cohomology) \cite{Kas, L2},
and they are called the global bases and canonical bases, respectively.
In \cite{GL}, it was shown that these two bases coincide with
each other.
Lusztig's construction was inspired by Ringel's realization of
$\U_v^{-}(\g)$ in terms of the Hall algebra of the underlying
Dynkin diagram.

\vspace{.1in}
The {\it generalized Kac-Moody algebras} were introduced by
Borcherds in his study of vertex algebras and
Monstrous Moonshine \cite{B88}.
In particular, the {\it Monster Lie algebra}, which is a special
example of generalized Kac-Moody algebras, played a crucial role
in the proof of the Moonshine conjecture \cite{B92}.
Recently, Jeong, Kang and Kashiwara developed the {\it crystal basis
theory} for quantum generalized Kac-Moody algebras, and constructed
the global bases for $\U_v(\g)$-modules in the category
$\mathcal{O}_{int}$ and for the subalgebra $\U_v^{-}(\g)$.

\vspace{.1in}

In this paper we provide the geometric construction of canonical
bases for quantum generalized Kac-Moody algebras.
We first extend Ringel's Hall algebra approach to the case of
quantum generalized Kac-Moody algebras.
In order to account for simple imaginary roots,
we consider the quivers in which edge loops are allowed,
and add some nilpotency relations at each such edge loop.
Next we generalize to this context Lusztig's construction of canonical bases
in terms of perverse sheaves on varieties of representations of quivers.

\vspace{.1in}

An important difference between our construction and the classical one
is that our canonical bases are made of \textit{semi}simple perverse
sheaves rather than simple perverse sheaves.
That this is necessary is already clear
in the case of a quiver consisting of a single loop,
where the corresponding generalized Kac-Moody algebra is the
Heisenberg algebra and the space of representations
is the nilpotent cone.
In particular, the canonical basis is not
(quasi)-orthonormal anymore but only (quasi)-orthogonal. This makes the extension
of Lusztig's construction more delicate.
We conjecture that the canonical bases thus obtained coincide
with the global bases constructed in \cite{JKK}.
We prove this conjecture in the case of generalized Kac-Moody algebras
with no isotropic simple roots, which includes many interesting
generalized Kac-Moody algebras such as the Monster Lie algebra.
We expect our work will lead to a variety of combinatorial and geometric
developments in the study of generalized Kac-Moody algebras
and their representations.

\vspace{.1in}

\noindent
\textbf{Acknowledgements.} \
Part of this work was done while both authors visited
Yale University, and while the second author visited
Korea Institute for Advanced Study. They thank both institutions
for their hospitality.
The first author is very grateful to Professor Young-Hoon Kiem
at Seoul National University for many helpful discussions.

\vspace{.2in}

\section{Quantum generalized Kac-Moody algebras}

\vspace{.1in}

\paragraph{\textbf{1.1. Generalized root datum.}}
Let $I$ be a countable index set.
In this paper, a \textit{generalized root datum} is a matrix
$A=(a_{ij})_{i,j \in I}$
satisfying the following conditions~:
\begin{enumerate}
\item[i)] $a_{ii} \in \{2,0,-2,-4,\ldots\}$,
\item[ii)] $a_{ij}=a_{ji} \in \Z_{\leq 0}$
\end{enumerate}

Such a matrix is a special case of \textit{Borcherds-Cartan} matrix (see
\cite{B88}).  Let $I^{re}=\{i \in I\;|a_{ii}=2\}$ and
$I^{im}=I\backslash I^{re}$.  We also assume that we are given a
collection of positive integers  (the \textit{charge} of $A$)
$\m=(m_i)_{i\in I}$ with $m_i=1$ whenever $i \in I^{re}$.

\vspace{.2in}

\paragraph{\textbf{1.2. Definition.}}
The {\it quantum generalized Kac-Moody algebra}
associated with $(A,\m)$ is the
(unital) $\C(v)$-algebra $\U_v(\g_{A,\m})$ (or simply $\U_v(\g)$)
 generated by the elements $K_i,K_i^{-1},
E_{i,k}, F_{i,k}$ for $i \in I$, $k=1, \ldots, m_i$ subject to the
following  set of relations :

\begin{equation}\label{E:01}
K_iK_i^{-1}=K_i^{-1}K_i=1, K_iK_j=K_jK_i,
\end{equation}
\begin{equation}\label{E:02}
K_iE_{jk}K_i^{-1}=v^{a_{ij}}E_{jk}, \qquad K_iF_{jkj}K_i^{-1}=
v^{-a_{ij}}F_{jk},
\end{equation}
\begin{equation}\label{E:03}
E_{ik}F_{jl} - F_{jl}E_{ik} = \delta_{lk}\delta_{ij}
\frac{K_i-K_i^{-1}}{v-v^{-1}},
\end{equation}
\begin{equation}\label{E:04}
\begin{split}
\sum_{n=0}^{1-a_{ij}}(-1)^n
\left[\begin{matrix} 1-a_{ij}\\n\end{matrix}\right] E_{ik}^{1-a_{ij}-n}
E_{jl}E^n_{ik}&=0 \qquad \forall\;i \in I^{re}, j \in I,\\
\sum_{n=0}^{1-a_{ij}}(-1)^n
\left[\begin{matrix} 1-a_{ij}\\n\end{matrix}\right] F_{ik}^{1-a_{ij}-n}
F_{jl}F^n_{ik}&=0 \qquad \forall\;i \in I^{re}, j \in I,
\end{split}
\end{equation}
\begin{equation}\label{E:05}
E_{ik}E_{jl} - E_{jl}E_{ik}=F_{ik}F_{jl} - F_{jl}F_{ik} =0 \qquad if\;
a_{ij}=0, \end{equation}
where as usual we put
$$[n]=\frac{v^n-v^{-n}}{v-v^{-1}}, \qquad [n]!=[2] \cdots [n], \qquad
\left[\begin{matrix} n\\k\end{matrix}\right]=\frac{[n]!}{[n-k]![k]!}.$$

\paragraph{}The algebra $\U_v(\g_{})$ is equipped with a Hopf algebra
structure as follows (see \cite{BKM, Ka95})~:
$$\Delta(K_i)=K_i \otimes K_i,$$
$$\Delta(E_{ik})=E_{ik} \otimes K_i^{-1} + 1 \otimes E_{ik},
\qquad \Delta(F_{ik})
=F_{ik} \otimes 1 + K_i \otimes F_{ik},$$
$$\epsilon(K_i)=1, \quad \epsilon(E_{ik})=\epsilon(F_{ik})=0,$$
$$S(K_i)=K_i^{-1}, \quad S(E_{ik})=-E_{ik}K_i, \quad S(F_{ik})=
-K_i^{-1}F_{ik}.$$

\vspace{.1in}

\paragraph{}Let  $\A=\Z[v,v^{-1}]$. The {\it integral form}
$\U_{\A}(\g_{})$ is the $\A$-subalgebra of $\U_v(\g_{})$ generated by
$K_j^{\pm 1}, E_{ik}, F_{ik}$ for $j \in I, i \in I^{im}$ and
$k=1, \ldots, m_i$, and the divided powers
$$E_{ik}^{(n)}=\frac{E_{ik}^n}{[n]!},\qquad
F_{ik}^{(n)}=\frac{F_{ik}^n}{[n]!},\qquad \text{for}\; i \in I^{re}\;
\text{and}\; n >0.$$
It follows from the definitions that the subalgebra $\U^{\geq
0}_v(\g_{})$ of $\U_v(\g_{})$ generated by $K_i^{\pm 1}$ and $F_{ik}$
for $i \in I$ and $k=1, \ldots, m_i$ is a Hopf subalgebra. We also
denote by $\U^-_v(\g_{})$ the subalgebra generated by $F_{ik}$ for $i
\in I$ and $k=1, \ldots, m_i$. We  will use similar notations for
$\U_{\A}(\g)$.

\vspace{.1in}

\paragraph{}Finally, let $u \mapsto \overline{u}$ be the semilinear
involution of $\U_v(\g)$ defined by $\overline{v}=v^{-1}$,
$\overline{F_{ik}}=F_{ik}$, $\overline{E_{ik}} =E_{ik}$ and
$\overline{K_i}=K_i^{-1}$.

\vspace{.2in}

\paragraph{\textbf{1.3. Nondegenerate bilinear form.}}
It follows from \cite{L1}
that there exists a unique symmetric bilinear form $\langle
\,,\,\rangle:  \U^{\leq 0}_v (\g) \otimes \U^{\leq 0}_v (\g) \to \C(v)$
satisfying $$\langle K_i, K_j \rangle =v^{-a_{ij}}, \quad \langle
F_{ik}, F_{jl}\rangle = \delta_{ij}\delta_{kl}, \quad \langle F_{ik},
K_{j} \rangle =0$$
and the invariance condition
$$\langle a, bc \rangle = \sum_n \langle a_n^{(1)}, b \rangle \langle
a_n^{(2)} , c\rangle.$$
Here, we have used Sweedler's notation $\Delta(a)=\sum_n a_n^{(1)}
\otimes a_n^{(2)}$.

\vspace{.1in}

\begin{prop}[\cite{SVdB}] The restriction of $\langle\,,\,\rangle$ to
$\U^-_v(\g)$ is nondegenerate.
\end{prop}

\vspace{.2in}

\section{Quivers and Hall algebras}

\vspace{.1in}

\paragraph{\textbf{2.1. Quivers.}} Let $Q$ be an arbitrary locally
finite quiver with vertex set $I$ and (oriented) edge set $\Omega$. For
$\sigma \in \Omega$
we denote by $o(\sigma)$ and $i(\sigma)$ the outgoing and incoming
edges,  respectively,
and sometimes use the notation $o(\sigma) \stackrel{\sigma}{\to}
i(\sigma)$.  We will denote by
$c_i$ the number of loops at $i$ (i.e., the number of edges $\sigma$
with  $i(\sigma)=o(\sigma)=~i$).

\vspace{.1in}

\paragraph{}Recall that a representation of $Q$ over a field $k$ is a
collection $(V_i, x_{\sigma})_{i \in I, \sigma \in \Omega}$, where $V_i$
is a finite-dimensional $k$-vector space, $V_i=0$ for almost all $i$ and
$x_{\sigma} \in \mathrm{Hom}_k (V_{o(\sigma)}, V_{i(\sigma)})$.  We
denote by $Rep_k(Q)$
the abelian category of representations of $Q$ over $k$. For any $M,N
\in  Rep_k(Q)$ the spaces $\mathrm{Hom}(M,N)$ and
$\mathrm{Ext}^1(M,N)$ are finite-dimensional
and we have $\mathrm{Ext}^n(M,N)=0$ for $n >1$.

\vspace{.1in}

Let us set $\mathbf{dim}(V_i,x_{\sigma})=
\mathbf{dim}(V_i)=(\mathrm{dim}_k V_i)
\in
\N^{\oplus I}$. The {\it Euler form} given by
$$\langle M, N \rangle =\mathrm{dim}_k
\mathrm{Hom}(M,N)-\mathrm{dim}_k
\mathrm{Ext}^1(M,N)$$
factors through $\Z^{\oplus I}$ and yields the
following nonsymmetric bilinear pairing
$\langle\,,\,\rangle : \Z^{\oplus I} \otimes_{\Z} \Z^{\oplus I} \to \Z$
defined by
$$((d_i),(d'_i)) \mapsto \sum_i (1-c_i)d_id_i' -\sum_{i\neq j}
r_{ij}d_id'_j,$$ where $r_{ij}$ denotes the number of edges going from
$i$ to $j$. We also  introduce the associated symmetric form
$$(a,b)=\langle a,b \rangle + \langle b,a\rangle.$$

\vspace{.1in}

\paragraph{}Observe that there is a correspondence between locally
finite quivers and generalized root datum : to $Q$ we associate the
matrix $A_Q =(a_{ij})_{i,j \in I}$ with $a_{ij}=(\epsilon_i,
\epsilon_j)$, where $(\epsilon_i)_i$ stands for the standard basis of
$\Z^{\oplus I}$.  In particular,
$a_{ii}=2-2c_i$, and $i \in I$ is real if and only if $c_i=0$.
If $V$ is an $I$-graded vector space, we will write $V=V^{re} \oplus
V^{im}$  for its decomposition according to $I=I^{re} \sqcup I^{im}$.
\vspace{.2in}

\paragraph{\textbf{2.2. Hall algebras.}} Let us now assume that $k$ is a
 finite field with $q$ elements and fix $v \in \C$ such that
$v^2=q^{-1}$. Let  us denote by $\mathcal{I}so_k(Q)$ the set of
isomorphism classes of objects in $Rep_k(Q)$, and let us write $[M] \in
\mathcal{I}so_k(Q)$ for the class of an object $M$. The {\it Hall
algebra} $\H_k(Q)$ is by definition the associative  $\C$-algebra with
basis indexed by $\mathcal{I}so_k(Q)$ and equipped with the  product

\begin{equation}\label{E:1}
[M]\cdot [N]=v^{-\langle \mathbf{dim}\;M,\mathbf{dim}\;N\rangle}
\sum_{[P] \in \mathcal{I}so_k(Q)} \mathcal{P}^P_{M,N}[P],
\end{equation}
where
$$\mathcal{P}^P_{M,N}=\#\{X \subset P\;|\;X \in [N], P/X \in [M]\}$$

\vspace{.1in}

Observe that the above quantity is finite since $\mathrm{Hom}(N,P)$ is
finite, and that
the sum (\ref{E:1}) only has finitely many terms since
$\mathrm{Ext}^1(M,N)$ is finite.

\vspace{.1in}

\paragraph{}Following Green \cite{Green}, we extend $\H_k(Q)$ by a
commutative  polynomial algebra : let $\dot{\H}_k(Q)$ be the
algebra generated by
$\H_k(Q)$ together with elements $\mathcal{K}_i^{\pm 1}$ for
$i \in I$ with the following relations :
\begin{equation}\label{E:2}
\mathcal{K}_i \mathcal{K}_i^{-1}=\mathcal{K}_i^{-1}\mathcal{K}_i=1,\quad
\mathcal{K}_i\mathcal{K}_j=\mathcal{K}_j\mathcal{K}_i,
\end{equation}
\begin{equation}\label{E:3}
\mathcal{K}_i[M]\mathcal{K}_i^{-1}=v^{-(\epsilon_i,
\mathbf{dim}\;M)}[M]. \end{equation}

It is easy to see that the multiplication map $\C[\mathcal{K}_i^{\pm 1}]
\otimes \H_k(Q) \to \dot{\H}_k(Q)$ is an isomorphism of vector spaces.
Finally, we equip the algebra $\dot{\H}_k(Q)$ with a coproduct by
setting
$\Delta(\mathcal{K}_i)=\mathcal{K}_i \otimes \mathcal{K}_i$ for all $i
\in I$ and
$$\Delta([M])=\sum_{N \subset M}
v^{-\langle \mathbf{dim}\;M/N, \mathbf{dim}\;N \rangle}
\mathcal{P}^M_{M/N,N} \frac{a_{M/N}a_N}{a_M}[M/N]
\mathcal{K}^{\mathbf{dim}\;N} \otimes [N],$$
where $a_M=\#Aut(M)$ and $\mathcal{K}^{\mathbf{dim}\;N}=\prod_i
\mathcal{K}_i^{(\mathbf{dim}\;N)_i}$.

\vspace{.1in}

\paragraph{}The following is one of the main results in \cite{Green} :

\vspace{.1in}

\begin{theo}[\cite{Green}] There exists a Hopf pairing
$\langle\,,\,\rangle_G$ on $\dot{\H}_k(Q)$ such that
$$\langle
\mathcal{K}_i,\mathcal{K}_j\rangle_G=v^{-(\epsilon_i,\epsilon_j)}, \quad
\langle [M],[N]\rangle_G=\delta_{[M],[N]}\frac{1}{a_M}.$$
\end{theo}

\vspace{.2in}

\paragraph{\textbf{2.3. Composition algebra.}} If $i \in I^{re}$, then
there exists a unique simple object $S_i \in Rep_k(Q)$ such that
$\mathbf{dim}\;S_i= \epsilon_i$. On the other hand, if $i \in I^{im}$
then the set of simple representations of dimension $\epsilon_i$ is in
bijection with $k^{c_i}$ : if $\sigma_1, \ldots, \sigma_{c_i}$ denote
the simple loops at $i$ then to $\underline{\lambda}=(\lambda_1, \ldots
,\lambda_{c_i})$ corresponds the simple module
$S_i(\underline{\lambda})=(V_j,x_{\sigma})$ with
$\mathrm{dim}_k V_j=
\delta_{ij}$ and $x_{\sigma_l}=\lambda_l Id$ for $l=1, \ldots, c_i$.

\vspace{.1in}

\paragraph{} For each $i \in I^{im}$, let us now assume that
we are given a positive integer $m_i$.
We choose $\underline{\lambda}_i^{(l)} \in k^{c_i}$ for $l=1, \ldots,
m_i$ in such a way that if $|k^{c_i}| \geq m_i$ then
$\underline{\lambda}_i^{(l)} \neq \underline{\lambda}^{(l')}_i$ for $l
\neq l'$. Consider the \textit{composition subalgebra} $\mathbf{C}_k(Q)
\subset \H_k(Q)$ generated by the classes $[S_i]$ for $i \in I^{re}$ and
$[S_{i,l}]:=[S_i(\underline{\lambda}_i^{(l)})]$ for $i \in I^{im}$ and
$l=1, \ldots, m_i$. We also define $\dot{\mathbf{C}}_k(Q)$ as the
subalgebra of  $\dot{\H}_k(Q)$
generated by $\mathbf{C}_k(Q)$ and $\C[\mathcal{K}_i^{\pm 1}]$.

\vspace{.1in}

\paragraph{}Following Ringel, we now define a {\it generic composition
algebra}: let $\mathbf{K}$ be an infinite set of (nonisomorphic) finite
fields, and let us choose for each $k \in \mathbf{K}$ an element $v_k
\in \C$  such that $v_k^2=|k|^{-1}$. Consider the direct product
$$\dot{\H}(Q)=\prod_{k \in \mathbf{K}} \dot{\H}_k (Q).$$

We view $\dot{\H}(Q)$ as a $\C[v,v^{-1}]$-module by mapping $v,v^{-1}$
to $(v_k)_k,(v_k^{-1})_k$ respectively.
We define $\dot{\mathbf{C}}(Q)$ to be the $\A$-subalgebra
of $\dot{\H}(Q)$ generated by $(\mathcal{K}_i)_k,
(\mathcal{K}_i^{-1})_k $, $([S_i])_k$ for $i \in I^{re}$ and
$([S_{i,l}])_k$ for $i \in I^{im}$ and $l=1, \ldots , m_i$. The
subalgebra  $\mathbf{C}(Q)$ is defined in a similar fashion.

\vspace{.1in}

Our first result is an extension of a well-known theorem of Green (in
the case of a quiver without loop) :

\vspace{.1in}

\begin{theo} The assignement $K_{i}^{\pm 1} \mapsto (\mathcal{K}_i^{\pm
1})_k$, $F_i \mapsto ([S_i])_k$ for $i \in I^{re}$ and $F_{il} \mapsto
([S_{i,l}])_k$ for $i \in I^{im}$ extends to an isomorphism of
bialgebras $\U_v^{\leq 0}(\g_{A_Q,\mathbf{m}}) \simeq
\dot{\mathbf{C}}(Q) \otimes \C(v)$  where $\mathbf{m}=(m_i)_i$.
\end{theo}

\vspace{.1in}

The proof of the theorem parallels the proof in \cite{Green} and is
given in  the next section.

\vspace{.2in}

\section{Hall algebra construction of $\U_{v}^{-}(\g)$}

\vspace{.1in}

\paragraph{}We first show that for any fixed field $k$, the assignement
$K_i^{\pm 1} \mapsto \mathcal{K}_i^{\pm 1}, F_i \mapsto [S_i], F_{il}
\mapsto [S_{i,l}]$ extends to a bialgebra homomorphism
$\varphi_k:\U_v^{\leq 0}(\g_{A_Q,\mathbf{m}}) \to
\dot{\mathbf{C}}_k(Q)$.
>From the relations (\ref{E:01}) and (\ref{E:02}), the equations
> (\ref{E:2}),
(\ref{E:3})
and the definitions of the coproduct, it is easy to see that it is
enough to check that we have a morphism $\varphi_k: \U^-_{v}(\g) \to
\mathbf{C}_k(Q)$. The corresponding computations are similar to the ones
in \cite{Rin}, which we reproduce for the reader's convenience.

\vspace{.1in}

\paragraph{}We will first check that if $i \in I^{re}$ and $j \in I$,
then $[S_i]$ and $[S_j]$ satisfy the $q$-Serre relation (\ref{E:04}). To
simplify the notations, let us set $r_1=r_{ij},r_2=r_{ji}$ and
$r=r_1+r_2$.

\vspace{.1in}

An easy induction shows that
$[S_i]^{(l)}:=\frac{[S_i]^l}{[l]!}=v^{-l(l-1)} [S_i^{\oplus l}]$. We
have
$$[S_i]^{(l)}[S_{j,t}]=v^{-l(l-1)+lr_1}\sum_{[T] \in
\mathcal{I}_1}[T],$$ where $\mathcal{I}_1=\{[T]\;|\;\exists X \subset T
\;s.t.\; X \simeq S_{j,t},  \; T/X \simeq S_i^{\oplus l}\}$. For a
representation $P$ of $Q$ of dimension $(r+1)\epsilon_i \oplus
\epsilon_j$, we define
$$U_P=\bigcap_{i \stackrel{\sigma}{\to} j} \mathrm{Ker}\;x_{\sigma},
\qquad  V_P=\sum_{j \stackrel{\sigma}{\to} i} \mathrm{Im}\;x_{\sigma},$$
and set $u_P=\mathrm{dim}\;U_P, v_P=\mathrm{dim}\;V_P$. For $a,b \in
\Z$,  let us denote by $Gr_a^b(k)$ the Grassmanian of $a$-dimensional
subspaces in  $k^b$ (the empty set if $a<0$ or $a>b$). A direct
computation now shows that $$[S_i]^{(l)}
[S_{j,t}][S_i]^{(n)}=v^{nr_2+lr_1-nl-l(l-1)-n(n-1)}\sum_{[P]}
\sigma_{P,n}[P],$$
where $\sigma_{P,n}=0$ unless $V_P \subset U_P$, in which case we have
$$\sigma_{P,n}=\#Gr_{n-v_P}^{u_P-v_P}=v^{-(u_P-n)(n-v_P)}
\left[\begin{matrix} u_P-v_P\\n-v_P\end{matrix}\right].$$

\vspace{.1in}

Setting $n=r+1-l$ and summing up, we obtain
$$\sum_{l=0}^{r+1} (-1)^l [S_i]^{(l)}[S_{j,t}][S_i]^{(n)}=
\sum_{[P]\;s.t.V_P \subset U_P} \gamma_P[P],$$
where
\begin{equation*}
\begin{split}
\gamma_P &=\sum_{l=0}^{k+1} (-1)^l
v^{nr_2+lr_1-nl-l(l-1)-n(n-1)-(u_P-n)(n-v_P)}
\left[\begin{matrix} u_P-v_P\\n-v_P\end{matrix}\right]\\
&=v^{-(r+1)r_2 +u_Pv_P}
\sum_{n=0}^{r+1}(-1)^{r+1-n}v^{(2r_2+1-u_P-v_P)n} \left[\begin{matrix}
u_P-v_P\\n-v_P\end{matrix}\right].
\end{split}
\end{equation*}

\vspace{.1in}

Observe that $u_P \geq r_2+1 >v_P$ for any $P$. We deduce that
$1-u_P-v_P \leq 2r_2+1-u_P-v_P \leq u_P +v_P -1$. The $q$-Serre relation
is now a  consequence of the following well-known identity
(see, for example, \cite[(3.2.8)]{Kas}).

\begin{lem} Let $m \geq 1$ and let $1-m \leq d \leq m-1$ with $d \equiv
m-1  \;(\mathrm{mod}\;2)$. Then
$$\sum_{n=0}^m (-1)^n v^{dn}
\left[\begin{matrix} m\\n\end{matrix}\right]=0.$$
\end{lem}

\vspace{.1in}

\paragraph{}Finally, let $i,j \in I$ such that $(\epsilon_i,
\epsilon_j)=0$. There are two possibilities: either $i \neq j$ and
$r_{ij}=r_{ji}=0$, or $i=j$ and $c_i=1$. In the first case, the relation
(\ref{E:05}) is obviously satisfied by $[S_{i,s}]$ and $[S_{j,t}]$. In
the second case, let us denote by $\sigma$ the edge loop at $i$. If
$\lambda_1 \neq \lambda_2$, then  the decomposition into
$x_{\sigma}$-eigenspaces shows that any short exact sequence $$0 \to
S_i(\lambda_1) \to M \to S_i(\lambda_2) \to 0$$
canonically splits; i.e.,
$[S_i(\lambda_1)][S_i(\lambda_2)]=[S_i(\lambda_1) \oplus
S_i(\lambda_2)]$, and the relation (\ref{E:05}) follows. Hence there is
a well-defined
algebra homomorphism $\varphi: \U^{\leq 0}_v(\g) \to \dot{\mathbf{C}}(Q)
\otimes \C(v)$.

\vspace{.1in}

\paragraph{} It remains to prove that the map $\varphi$ is injective.
For this,  recall the following lemma of Green (\cite{Green}). Set
$\overline{I}= \{(i,l)\;|\; 1 \leq l \leq m_i\}$. For $\nu \in
\N^{\overline{I}}$, let $I_{\nu}$ be the set of sequences
$\underline{a}=(a_1, \ldots, a_n)$ of  elements of $\overline{I}$ such
that $\#\{h\;|a_h=(i,l)\}=\nu_{(i,l)}$ for every $(i,l) \in
\overline{I}$.

\vspace{.1in}

\begin{lem}[\cite{Green}] Let $\nu \in \N^{\overline{I}}$, and let
$\underline{a}, \underline{b} \in \mathcal{I}_{\nu}$. There exists a
polynomial $M_{\underline{a},\underline{b}}(t) \in \Z[t,t^{-1}]$ such
that
$$\langle F_{a_1}\cdots F_{a_n}, F_{b_1}, \cdots F_{b_m}\rangle=
M_{\underline{a},\underline{b}}(v)$$
and such that for every finite field $k \in \mathbf{K}$,

$$\langle [S_{a_1}] \cdots [S_{a_n}], [S_{b_1}] \cdots [S_{b_m}]
\rangle_G =M_{\underline{a},\underline{b}}(v_k) \prod_{(i,l) \in
\overline{I}} \langle [S_{i,l}], [S_{i,l}] \rangle_G^{\nu_{(i,l)}}.$$
\end{lem}

\vspace{.1in}

\paragraph{} Now let $x=\sum_{\underline{a}} c_{\underline{a}}(v)F_{a_1}
\cdots F_{a_n} \in \mathrm{Ker}\;\varphi$, so that for any $k \in
\mathbf{K}$ we have

$$\varphi_k(x)=\sum_{\underline{a}}c_{\underline{a}}(v_k) [S_{a_1}]
\cdots [S_{a_n}]=0.$$
In particular, for any $\underline{b}$, we have
$$\sum_{\underline{a}}
c_{\underline{a}}(v_k)M_{\underline{a},\underline{b}} (v_k)=0.$$
Since $\mathbf{K}$ is infinite, this implies that
$$\sum_{\underline{a}} c_{\underline{a}}(v)
M_{\underline{a},\underline{b}}(v) =0 \in \C[v,v^{-1}].$$
Thus $x$ lies in the radical of $\langle\;,\;\rangle$. But by
Proposition 1.1., the form $\langle\;,\;\rangle$ is nondegenerate. Thus
$x=0$ and Theorem 2.2 is proved. \qed

\vspace{.2in}

\section{The Algebra of Semisimple Perverse Sheaves}

\vspace{.1in}

\paragraph{\textbf{4.1. Quiver representation varieties.}} We keep the
notations of Section 2.1., but assume that $k=\overline{\mathbb{F}_q}$.
For simplicity, we will only consider the
case of a generalized Kac-Moody algebra $\g_{A_Q,\mathbf{m}}$ with
trivial charge; i.e., $m_i=1$ for all $i \in I$ (in fact, this is not
restrictive; see \cite[Remark 1.3]{JKK}).

\vspace{.1in}

For all $\mathbf{d} \in \N^{\oplus I}$, we fix an $I$-graded $k$-vector
spaces $V_{\d}=\bigoplus_iV_i$ such that
$\mathbf{dim}(V_\d)=\mathbf{d}$. Let  $\mathbf{d}^{im}$ and
$\mathbf{d}^{re}$ be the imaginary and real components  of $\mathbf{d}$,
respectively. Denote by
$$E_{\mathbf{d}}=\big\{(x_{\sigma})\;|x_{\sigma_1} \cdots
x_{\sigma_N}=0\; \mathrm{for\;any\;}\sigma_i \in \Omega\;\mathrm{and}\;
N \gg 0\big\} \subset \bigoplus_{\sigma \in \Omega}
\mathrm{Hom}\;(V_{o(\sigma)}, V_{i(\sigma)})$$
the set of \textit{nilpotent} representations of $Q$ in $V_\d$.
Now let $\mathbf{i}=(i_1, \ldots, i_r)$ be a sequence of vertices
$i_l \in I$ such that $\sum \epsilon_{i_l}=\mathbf{d}$.
Consider the variety of $I$-graded flags
$$\mathcal{F}_{\mathbf{i}}=\{D_{\bullet}\;|0 \subseteq D_1 \subseteq
\cdots \subseteq D_r=V_\d;\; \mathbf{dim}
(D_l/D_{l-1})=\epsilon_{i_l}\}.$$ Finally, we define the {\it incidence
varieties}
$$\widetilde{\mathcal{F}}_{\mathbf{i}}=\{(x,D_{\bullet})\;|x(D_i)
\subseteq D_{i-1}\} \subset E_{\mathbf{d}} \times
\mathcal{F}_{\mathbf{i}},$$
$$\widetilde{\mathcal{F}}^{im}_{\mathbf{i}}=\{(x,D_{\bullet})\;|x(D_i)
\subseteq D_{i-1} \oplus \bigoplus_{i \in I^{re}} V_i\}
\subset E_{\mathbf{d}} \times \mathcal{F}_{\mathbf{i}}^{im},$$
where $\mathcal{F}_{\mathbf{i}}^{im}$ is defined
as $\mathcal{F}_{\mathbf{i}}$ by replacing $\mathbf{d}$ by
$\mathbf{d}^{im}$.

\vspace{.1in}

Thus we have a commutative diagram, with obviously defined projection
maps\,: $$
\xymatrix{
\widetilde{\mathcal{F}}_{\mathbf{i}}
\ar[r]^{\pi_1'} \ar[d]_{\pi'_2}&
\widetilde{\mathcal{F}}_{\mathbf{i}}^{im}  \ar[d]_{\pi_2}
\ar[r]^{\pi_1}&
E_{\mathbf{d}} \\
\mathcal{F}_{\mathbf{i}} \ar[r] &
\mathcal{F}_{\mathbf{i}}^{im}
}
$$
Note that $\pi_1$ and $\pi_1'$ are smooth proper,
while $\pi_2$ and $\pi_2'$ are vector bundles.

\vspace{.2in}

\paragraph{\textbf{4.2. Notations.}} We use the notations in
\cite[Chapter 8]{L1} regarding perverse sheaves. In particular, for an
algebraic variety $X$ defined over $k$ we denote by $\mathcal{D}(X)$
(resp. $\mathcal{Q}(X)$, resp. $\mathcal{M}(X)$) the derived category of

$\overline{\mathbb{Q}}_l$-constructible sheaves (resp. the category of
semisimple $\overline{\mathbb{Q}}_l$-constructible complexes, resp. the
category of perverse sheaves) on $X$.
The Verdier dual of a complex $P$ is
denoted by $D(P)$ and perverse cohomology by $H^{\bullet}(P)$.
If $G$ is a connected algebraic group
acting on $X$, then we denote by
$\mathcal{Q}_G(X)$ and $\mathcal{M}_G(X)$
the corresponding categories of $G$-equivariant complexes.
Finally, for $P_1, P_2 \in \mathcal{D}_{G}(X)$ and $j \in \Z$ we set
$D_j(P_1,P_2)$ to be the dimension of the space denoted
$\mathbf{D}_j(X,G,P_1,P_2)$ in \cite{L1}.

\vspace{.2in}

\paragraph{\textbf{4.3. Induction and restriction functors.}} The group
$G_{\mathbf{d}}=\prod_i GL(V_i)$ naturally acts on $E_{\mathbf{d}}$.
Following Lusztig (see \cite{L1}), let us fix an embedding of $I$-graded
vector spaces $V_{\d_1} \to V_{\d_1+\d_2}$ and an isomorphism
$i:V_{\d_1+\d_2} /V_{\d_1} \simeq V_{\d_2}$, and consider the diagram
$$\xymatrix{
 E_{\d_1+\d_2} & F \ar[r]^-{\kappa} \ar[l]_-{\iota} & E_{\d_1} \times
E_{\d_2}}, $$
where $F=\{x \in E_{\d_1+\d_2}\;| x(V_{\d_1}) \subset V_{\d_1}\}$,
$\iota$ is  the canonical embedding and
$\kappa(x)=(x_{|V_{\d_1}}, i_*(x_{|V_{\d_1+\d_2}/V_{\d_1}}))$. It is
clear that $\kappa$ is a vector bundle.

Similarly, consider the diagram
$$\xymatrix{
E_{\d_1} \times E_{\d_2} & E' \ar[l]_-{p_1} \ar[r]^-{p_2} & E''
\ar[r]^-{p_3}& E_{\d_1 + \d_2}},$$
where
\begin{enumerate}
\item[-]
$E'$ is the variety of all quadruples $(x,V,\alpha,\beta)$ satisfying $x
\in E_{\d_1+\d_2}; \;V \subset V_{\d_1 + \d_2};\;
x(V) \subset V;\; \alpha: V \simeq V_{\d_1};\; \beta: V_{\d_1+\d_2} / V
\simeq V_{\d_2}$,
\item[-]
$E''$ is the variety of pairs $(x,V)$ satisfying $x \in E_{\d_1+\d_2};\;
 V \subset V_{\d_1 + \d_2};\;
x(V) \subset V;\; \mathbf{dim}(V) = \d_1,$
\item[-]
$p_1(x,V,\alpha,\beta)=(\alpha_* (x_{|V}),
\beta_*(x_{|V_{\d_1+\d_2}/V}))$,
\item[-]
$p_2(x,V,\alpha,\beta)=(x,V)$ and $p_3(x,V)=x$.
\end{enumerate}
Note that $G_{\d_1} \times G_{\d_2}$ acts on $E_{\d_1}\times E_{\d_2}$
and on $E'$; that $G_{\d_1+\d_2}$ acts on $E', E'', E_{\d_1+\d_2}$ and
trivially on $E_{\d_1} \times E_{\d_2}$; that all maps are equivariant
for these groups. Also observe that $p_1$ is smooth with connected
fibers, $p_2$ is a principal  $G_{\d_1} \times G_{\d_2}$-bundle while
$p_3$ is proper.

\vspace{.1in}

Define the functor
$$\widetilde{\mathrm{Res}}_{\d_1,\d_2}^{\d_1+\d_2}=\kappa_!\iota^*:
\mathcal{Q}_{G_{\d_1+\d_2}}
(E_{\d_1+\d_2}) \to \mathcal{D}(E_{\mathbf{d}_1}\times
E_{\mathbf{d}_2}),$$ and put
$\mathrm{Res}_{\mathbf{d}_1+\mathbf{d}_2}^{\mathbf{d}_1 ,
\mathbf{d}_2}=
\widetilde{\mathrm{Res}}_{\mathbf{d}_1,\mathbf{d}_2}^{\mathbf{d}_1+
\mathbf{d}_2}[l_1-l_2-2\sum_i (\d_1)_i(\d_2)_i]$, where
$l_1$ and $l_2$ are the dimensions of the
fibres of $p_1$ and $p_2$, respectively.
Also, we define the functor (see \cite[Section 9.2]{L1}) :
$$\widetilde{\mathrm{Ind}}_{\mathbf{d}_1,\mathbf{d}_2}^{\mathbf{d}_1 +
\mathbf{d}_2}
= p_{3!}p_{2\flat}p_1^* :
\mathcal{Q}_{G_{\d_1} \times G_{\d_2}}(E_{\mathbf{d}_1}\times
E_{\mathbf{d}_2}) \to \mathcal{D}(E_{\d_1+\d_2}),$$
and set $\mathrm{Ind}_{\mathbf{d}_1,\mathbf{d}_2}^{\mathbf{d}_1 +
\mathbf{d}_2}=
\widetilde{\mathrm{Ind}}_{\mathbf{d}_1,\mathbf{d}_2}^{\mathbf{d}_1 +
\mathbf{d}_2}[l_1-l_2]$.
With this shift, the functor
$\mathrm{Ind}$ commutes with Verdier duality (\cite[\S 9.2.5]{L1}).

\vspace{.2in}

\paragraph{\textbf{4.4. A class of semisimple perverse sheaves.}}
Let $\mathbf{i}$ be a sequence of vertices as in Section~4.1.
The variety $\widetilde{\mathcal{F}}_{\mathbf{i}}$ being smooth, the
constant sheaf
$(\overline{\mathbb{Q}}_l)_{\widetilde{\mathcal{F}}_{\mathbf{i}}}
[\mathrm{dim}\;\widetilde{\mathcal{F}}_{\mathbf{i}}]$ is perverse. The
map $\pi_1': \widetilde{\mathcal{F}}_{\mathbf{i}} \to
\widetilde{\mathcal{F}}^{im}_{\mathbf{i}}$ is proper and
$G_{\d}$-equivariant, so by \cite{BBD}
the complex $\mathcal{L}'_{\mathbf{i}}=
\pi_{1!}'((\overline{\mathbb{Q}}_l)_{\widetilde{\mathcal{F}}_{\mathbf{i}}})
[\mathrm{dim}\;\widetilde{\mathcal{F}}_{\mathbf{i}}]$
is semisimple, $G_\d$-equivariant and satisfies
$D(\mathcal{L}'_{\mathbf{i}})= \mathcal{L}'_{\mathbf{i}}$. Let
$\mathcal{T}_{\mathbf{i}}$ be the set of all simple perverse sheaves
appearing (possibly with a shift) in
$\mathcal{L}'_{\mathbf{i}}$.

\vspace{.1in}

\begin{prop} If $P \in \mathcal{T}_{\i}$, then $\pi_{1!}(P)$
is a semisimple perverse sheaf. Moreover, if $c_i >1$ for all $i \in
I^{im}$, then $\pi_{1!}(P)$ is simple.\end{prop}
\noindent
\textit{Proof.} The semisimplicity of $\pi_{1!}(P)$ follows from the
Decomposition Theorem in \cite{BBD} and the fact that $\pi_1 \circ
\pi_1'$ is proper. Let us prove that $\pi_{1!}(P)$ is in addition
perverse.
Let $E'_{\d}\subset\bigoplus_{i \neq j} \mathrm{Hom}(V_i,V_j)^{r_{ij}}$
and $E_{\d}''\subset\bigoplus_{i} \mathrm{Hom}(V_i,V_i)^{c_i}$ be the
set of nilpotent representations and let
$u: E_{\d} \to E'_{\d}$, $t: E_{\d} \to E''_{\d}$ be the projections.
Finally, we set $\mathcal{G}_{\i}= (u \times
Id)(\widetilde{\mathcal{F}}_{\i}) \subset E'_{\d} \times
\mathcal{F}_{\i}$ and
$\mathcal{G}^{im}_{\i}= (u \times Id)(\widetilde{\mathcal{F}}^{im}_{\i})
\subset E'_{\d} \times \mathcal{F}^{im}_{\i}$, so that we have a
commutative  diagram
$$\xymatrix{
\widetilde{\mathcal{F}}_{\mathbf{i}}
\ar[r]^{\pi_1'} \ar[d]_{u \times Id}&
\widetilde{\mathcal{F}}_{\mathbf{i}}^{im}
\ar[d]^{u\times Id} \ar[r]^{\pi_1}&
E_{\mathbf{d}} \\
\mathcal{G}_{\mathbf{i}} \ar[r]^{s} &
\mathcal{G}_{\mathbf{i}}^{im}
}$$
Observe that the vertical maps are vector bundles, and that
$u \times Id:\widetilde{\mathcal{F}}_{\mathbf{i}} \to \mathcal{G}_{\i}$
is the pullback by $s$ of the bundle $u\times Id:
\widetilde{\mathcal{F}}_{\mathbf{i}}^{im}
\to \mathcal{G}_{\i}^{im}$. Hence
$\pi'_{1!}((\overline{\mathbb{Q}}_l)_{\widetilde{\mathcal{F}}_{\i}})=
\pi_{1!}'(u\times Id)^*((\overline{\mathbb{Q}}_l)_{\mathcal{G}_{\i}})=
(u\times Id)^*s_!((\overline{\mathbb{Q}}_l)_{\mathcal{G}_{\i}})$.  In
particular, any
of the simple perverse sheaves in $\mathcal{T}_\i$ is of the form
$IC(X,\mathfrak{L})$ with $X = (u\times Id)^{-1}(Y)$ for a smooth
irreducible  subvariety $Y \subset \mathcal{G}_{\i}^{im}$ and
$\mathfrak{L}=
(u \times Id)^* \mathfrak{K}$ for an irreducible local system
$\mathfrak{K}$ on $Y$. Let us show that
the restriction of $\pi_1$ to $\overline{X}$ is {\it semismall}; i.e.,
$\mathrm{dim}\;\overline{X} \underset{E_{\d}}{\times}
\overline{X}=\mathrm{dim} \; \overline{X}$.
Then standard arguments would show that $\pi_{1!}(IC(X,\mathfrak{L}))$
is perverse (see \cite {BM}). \\

\paragraph{} Let us denote by $\mathcal{O}_{\lambda} \subset
\mathfrak{sl}_{n}$ the nilpotent orbit associated to a partition
$\lambda$ of an integer $n$, and set
$s_{\lambda}=\mathrm{dim}\;\{D_{\bullet} \in
\widetilde{\mathcal{F}}_{n}\;| x(D_i) \subset D_{i-1}\}$ to be the
dimension of the Springer fiber over any point $x \in
\mathcal{O}_{\lambda}$. We say that a nilpotent element
$x_i=(x_i^j)_{j=1}^{c_i} \in \mathrm{Hom}(V_i,V_i)^{c_i}$ is of type
$\underline{\lambda}_i=(\lambda_i^1, \cdots, \lambda_i^{c_i})$ if $x_i^j
\in  \mathcal{O}_{\lambda_i^j}$.
Finally, we say that $x=(x_i)_i \in \bigoplus_i \mathrm{Hom}
(V_i,V_i)^{c_i}$ is of type
$\underline{\lambda}=(\underline{\lambda}_i)_i$ if $x_i$ is of type
$\underline{\lambda}_i$ for all $i$.\\

\paragraph{}The space $\overline{X}$ admits a finite partition
$\overline{X}= \bigsqcup_{\underline{\lambda}}
\overline{X}_{\underline{\lambda}}$, where
$$\overline{X}_{\underline{\lambda}}=\{(x, D_{\bullet}) \in
\overline{X}\;|\;t(x)\;\text{is\;of\; type\;} \underline{\lambda}\}.$$
Thus $\overline{X} \underset{E_{\d}}{\times}
\overline{X}=\bigsqcup_{\underline{\lambda}}
\overline{X}_{\underline{\lambda}} \underset{E_{\d}}{\times}
\overline{X}_{\underline{\lambda}}$. For each
$\underline{\lambda}=(\underline{\lambda}_i)_i$ let us fix partitions
$\lambda^{(i)}$ such that
\begin{enumerate}
\item[i)] $\mathrm{dim}\;\mathcal{O}_{\lambda_i^j} \leq \mathrm{dim}\;
\mathcal{O}_{\lambda^{(i)}}$ for all $j=1, \ldots, c_i$,
\item[ii)] $\lambda^{(i)}=\lambda_i^j$ for at least one $j \in \{1,
\ldots, c_i\}$. \end{enumerate}
It is clear that
\begin{equation}\label{E:canrec1}
\mathrm{dim}\;\overline{X}_{\underline{\lambda}} \underset{E_\d}{\times}
 \overline{X}_{\underline{\lambda}}
\leq \mathrm{dim}\;\overline{X}_{\underline{\lambda}} + \sum_i
s_{\lambda^{(i)}}. \end{equation}
On the other hand, since $\overline{X}=(u \times Id)^{-1}
(\overline{Y})$, we have \begin{equation*}
\mathrm{codim}_{\overline{X}}\; \overline{X}_{\underline{\lambda}} \geq
\sum_{i,j}  \mathrm{codim}_{\mathfrak{n}_i}  (\mathcal{O}_{\lambda_i^j}
\cap \mathfrak{n}_i),
\end{equation*}
where $\mathfrak{n}_i \subset \mathfrak{sl}(V_i)$ is the nilpotent
radical of any  Borel subalgebra. Moreover, by a well-known result of
\cite{Sp},
$$\mathrm{codim}_{\mathfrak{n}_i}(\mathcal{O}_{\lambda_i^j} \cap
\mathfrak{n}_i)= \mathrm{dim}\;\mathfrak{n}_i - \frac{1}{2}
\mathrm{dim}\;\mathcal{O}_{\lambda_i^j} = s_{\lambda_i^j},$$
which yields
\begin{equation}\label{E:canrec2}
\begin{split}
\mathrm{codim}_{\overline{X}}\; \overline{X}_{\underline{\lambda}} &
\geq \sum_{i,j} s_{\lambda_i^j}\\
&\geq \sum_i c_i s_{\lambda^{(i)}}.
\end{split}
\end{equation}

 Thus, combining (\ref{E:canrec1}) with (\ref{E:canrec2}) we obtain
\begin{equation}
\begin{split}
\mathrm{dim}\;\overline{X} \underset{E_{\d}}{\times} \overline{X} =&\;
\mathrm{sup}_{\underline{\lambda}} \;\mathrm{dim}\;
\overline{X}_{\underline{\lambda}} \underset{E_{\d}}{\times}
\overline{X}_{\underline{\lambda}}\\
\leq&\; \mathrm{sup}_{\underline{\lambda}}\; \mathrm{dim}\;
\overline{X}_{\underline{\lambda}} + \sum_i s_{\lambda^{(i)}}\\
\leq& \;\mathrm{sup}_{\underline{\lambda}}\; \mathrm{dim}\;
\overline{X}_{\underline{\lambda}} + \sum_i (1-c_i)s_{\lambda^{(i)}}\\
\leq& \;\mathrm{dim}\;\overline{X}
\end{split}
\end{equation}
as desired. Now assume that $c_i >1$ for all $i \in I^{im}$. Then
$\mathrm{dim}\; \overline{X}_{\underline{\lambda}}
\underset{E_{\d}}{\times} \overline{X}_{\underline{\lambda}} <
\mathrm{dim}\;\overline{X}$ as soon as $\O_{\lambda^{(i)}}$ is not
regular for some $i$. On the other hand, if $\O_{\lambda^{(i)}}$ is
regular for all $i$, then we have $ \overline{X}_{\underline{\lambda}}
\underset{E_{\d}}{\times} \overline{X}_{\underline{\lambda}} \subset
\Delta_{\overline{X}}$, where $\Delta_{\overline{X}}$ is the diagonal of
$\overline{X}$. Hence $\pi_1$ is  {\it small} (see \cite {BM}) and thus
$\pi_{1!}(P)$ is simple.\qed

\vspace{.2in}

\paragraph{\textbf{Definition.}} Set $\mathcal{P}_{\mathbf{i}}=
\{\pi_{1!}(P)\;|P \in \mathcal{T}_\mathbf{i}\}$. By the previous Lemma,
$\mathcal{P}_{\mathbf{i}}$ consists of semisimple
$G_{\d}$-equivariant perverse sheaves. Let $\mathcal{P}_{\d}=
\bigcup_{\mathbf{i}} \mathcal{P}_{\mathbf{i}}$ where the sum ranges over
all  sequences $\mathbf{i}$ such that $\sum_l \epsilon_{i_l} =\d$.
Further,   denote by $\mathcal{Q}_{\d}$ the category of complexes which
are
sums of shifts of elements in $\mathcal{P}_{\d}$. Thus
$\mathcal{Q}_{\d}$ is a full subcategory of
$\mathcal{Q}_{G_{\d}}(E_{\d})$. Finally, for $\d_1, \d_2 \in \N^{\oplus
I}$ we let $\mathcal{Q}_{\d_1}  \boxtimes \mathcal{Q}_{\d_2}$ be the
full subcategory of
$\mathcal{Q}_{G_{\d_1} \times G_{\d_2}}(E_{\d_1} \times E_{\d_2})$
consisting of sums of complexes $P_1 \boxtimes P_2$ with $P_1 \in
\mathcal{Q}_{\d_1}, P_2 \in \mathcal{Q}_{\d_2}$.

\vspace{.1in}

\begin{lem} The functors $\mathrm{Ind}$ and $\mathrm{Res}$ restrict to
functors
$$\mathrm{Ind}_{\d_1,\d_2}^{\d_1+\d_2} : \mathcal{Q}_{\d_1} \boxtimes
\mathcal{Q}_{\d_2} \to \mathcal{Q}_{\d_1+\d_2},$$
$$\mathrm{Res}_{\d_1,\d_2}^{\d_1+\d_2} : \mathcal{Q}_{\d_1+\d_2} \to
 \mathcal{Q}_{\d_1} \boxtimes
\mathcal{Q}_{\d_2}.$$
\end{lem}

\noindent
\textit{Proof.} To prove the first statement, it is enough to show that
 for any $P_1 \in
\mathcal{T}_{\mathbf{i}_1}, P_2 \in \mathcal{T}_{\i_2}$
we have $\mathrm{Ind}_{\d_1,\d_2}^{\d_1+\d_2} ( \pi_{1!} P_1 \boxtimes
\pi_{1!} P_2) \in \mathcal{Q}_{\d_1+\d_2}$.
For this, consider the following diagram
$$\xymatrix{
\widetilde{\mathcal{F}}_{\i_1}^{im} \times
\widetilde{\mathcal{F}}_{\i_2}^{im}  & E_{im}' \ar[l]_-{\rho_1}
\ar[r]^-{\rho_2}
 & E_{im}'' \ar[r]^-{\rho_3}&
\widetilde{\mathcal{F}}_{\i_1\i_2}^{im}
}$$
where $\i_1\i_2$ is the concatenation of the sequences $\i_1$ and
$\i_2$, and \begin{enumerate}
\item[-] $E'_{im}$ is the variety of quintuples $(x,V,D_{\bullet},
\alpha,\beta)$ such that $x \in E_{\d_1+\d_2};\; V \subset
V_{\d_1+\d_2};\; x(V) \subset V;\; (x,D_{\bullet}) \in
\widetilde{\mathcal{F}}_{\i_1\i_2}^{im};\;
V^{im} \in D_{\bullet};\; \alpha: V \simeq
V_{\d_1};\; \beta: V_{\d_1+\d_2}/V \simeq V_{\d_2}$,
\item[-] $E''_{im}$ is the variety of triples $(x,V,D_{\bullet})$
satisfying $x \in E_{\d_1+\d_2};\; V \subset V_{\d_1+\d_2};\;$ $x(V)
\subset V;\; (x,D_{\bullet}) \in
\widetilde{\mathcal{F}}_{\i_1\i_2}^{im};\;
V^{im} \in D_{\bullet};\;\mathbf{dim}\;V=\d_1$,
\item[-] $\rho_1(x,V,D_{\bullet},\alpha,\beta)=(\alpha_{*}(x_{|V},
D_{\bullet} \cap V),\beta_{*}(x_{|V_{\d_1+\d_2}/V},D_{\bullet}/V))$,
\item[-] $\rho_2(x,V,D_{\bullet},\alpha,\beta)=(x,V,D_{\bullet})$ and
$\rho_3(x,V,D_{\bullet})=(x,D_{\bullet})$.
\end{enumerate}
We set
$$\widetilde{\mathrm{Ind}}_{\i_1,\i_2}^{\i_1\i_2}=
\rho_{3!}\rho_{2\flat}\rho_{1}^*:
\mathcal{Q}_{G_{\d_1} \times
G_{\d_2}}(\widetilde{\mathcal{F}}_{\i_1}^{im} \times
\widetilde{\mathcal{F}}_{\i_2}^{im}) \to \mathcal{Q}_{G_{\d_1+\d_2}}
(\widetilde{\mathcal{F}}_{\i_1\i_2}^{im}).$$

\vspace{.1in}

\noindent
\textbf{Claim.} We have $\pi_{1!}
\circ (\widetilde{\mathrm{Ind}}^{\i_1\i_2}_{\i_1,\i_2})
=\widetilde{\mathrm{Ind}}_{\d_1,\d_2}^{\d_1+\d_2}\circ
(\pi_{1!} \boxtimes \pi_{1!})$.

\vspace{.1in}

\noindent
\textit{Proof of claim.} We have a commutative diagram
$$\xymatrix{
\widetilde{\mathcal{F}}_{\i_1}^{im} \times
\widetilde{\mathcal{F}}_{\i_2}^{im} \ar[d]^-{\pi_1 \times \pi_1}
& E_{im}' \ar[l]_-{\rho_1} \ar[r]^-{\rho_2} \ar[d]^-{\rho_4}
 & E_{im}'' \ar[r]^-{\rho_3} \ar[d]^-{\rho_5}&
\widetilde{\mathcal{F}}_{\i_1\i_2}^{im}
\ar[d]^-{\pi_1}\\
E_{\d_1} \times E_{\d_2} & E' \ar[l]_-{p_1} \ar[r]^-{p_2} & E''
\ar[r]^-{p_3}& E_{\d_1 + \d_2}
}$$
where $\rho_4$ and $\rho_5$ are the obvious projections. The leftmost
commutative square is a cartesian (pull-back) diagram, and therefore
$\rho_{4!}\rho_1^*=p_1^*(\pi_{1!}\boxtimes \pi_{1!})$. Similarly,
we have $p_2^*\rho_{5!}=\rho_{4!}\rho_2^*$, and hence
$p_{2\flat}\rho_{4!}=\rho_{5!} \rho_{2\flat}$. Thus we get
\begin{equation*}
\begin{split}
\widetilde{\mathrm{Ind}}_{\d_1,\d_2}^{\d_1+\d_2}(\pi_{1!}(P_1)
\boxtimes \pi_{1!}(P_2))&=p_{3!}p_{2\flat}p_1^*
(\pi_{1!}(P_1) \boxtimes \pi_{1!}(P_2))\\
&=p_{3!}p_{2\flat}
\rho_{4!}(\rho_1^*(P_1 \boxtimes P_2))\\
&=p_{3!}\rho_{5!}\rho_{2\flat}
(\rho_1^*(P_1 \boxtimes P_2))\\
&=\pi_{1!}\rho_{3!}\rho_{2\flat}
(\rho_1^*(P_1 \boxtimes P_2))\\
&=\pi_{1!}\big(
\widetilde{\mathrm{Ind}}_{\i_1,\i_2}^{\i_1\i_2}(P_1 \boxtimes P_2)\big).
\end{split}
\end{equation*}
This proves the claim.

\vspace{.1in}

Now, as in \cite[Lemma 9.2.3]{L1}, one can prove that
$\widetilde{\mathrm{Ind}}_{\i_1,\i_2}^{\i_1\i_2}(\mathcal{L}'_{\i_1}
\boxtimes \mathcal{L}'_{\i_2})=\mathcal{L}'_{\i_1\i_2}$. Thus
$\widetilde{\mathrm{Ind}}_{\i_1,\i_2}^{\i_1\i_2}(P_1 \boxtimes P_2)$ is
a sum of shifts of simple perverse sheaves appearing in
$\mathcal{L}'_{\i_1\i_2}$. It follows that
$\pi_{1!}\mathrm{Ind}_{\i_1,\i_2}^{\i_1\i_2}(P_1
\boxtimes P_2) \in \mathcal{Q}_{\d_1+\d_2}$. This proves the first part
of Lemma~4.2.

\vspace{.1in}

\paragraph{}The second statement is proved in a similar way: we consider
the diagram
$$\xymatrix{
 \widetilde{\mathcal{F}}_{\i} \ar[d]^-{\pi_1'} & F'' \ar[l]_-{\iota''}
\ar[d]^-{r} &\\
\widetilde{\mathcal{F}}_{\i}^{im} \ar[d]^-{\pi_1} & F' \ar[l]_-{\iota'}
\ar[d]^-{q}& \\
 E_{\d_1+\d_2} & F \ar[r]^-{\kappa} \ar[l]_-{\iota} &E_{\d_1} \times
E_{\d_2}} $$
where $F'$ (resp. $F''$) is the variety of pairs $(x,D_{\bullet}) \in
\widetilde{\mathcal{F}}_{\i}^{im}$ (resp. $(x,D_{\bullet}) \in
\widetilde{\mathcal{F}}_{\i}$) such that $x(V_{\d_1}) \subset
V_{\d_1}$.
Let $P \in \mathcal{T}_{\i}$. We have
$\widetilde{\mathrm{Res}}_{\d_1,\d_2}^{\d_1+\d_2}(\pi_{1!}P)=\kappa_!\iota^*
\pi_{1!} P =\kappa_!q_! (\iota')^* P$. Note that by \cite[Section
9.2]{L1},
$\widetilde{\mathrm{Res}}_{\d_1,\d_2}^{\d_1+\d_2}(\pi_{1!}\mathcal{L}'_{\i})$
is semisimple, hence so is $\kappa_!q_! (\iota')^* P$.\\

For a pair of complementary subsequences $\j_1$ and $\j_2$in
$\j=\i^{im}$, set $$F'(\j_1,\j_2)=\{(x,D_{\bullet}) \in F'\;|
D_{\bullet} \cap V_{\d_1} \in \mathcal{F}_{\j_1}, D_{\bullet}/
D_{\bullet}\cap V_{\d_1} \in
\mathcal{F}_{\j_2}\},$$
and
\begin{align*}
\kappa_{\j_1,\j_2}: F'(\j_1,\j_2) &\to
\widetilde{\mathcal{F}}_{\j_1}^{im} \times
\widetilde{\mathcal{F}}_{\j_2}^{im}\\
(x,D_{\bullet}) &\mapsto \big( (x_{|V_{\d_1}}, D_{\bullet} \cap
V_{\d_1}), (x_{|V_{\d}/V_{\d_1}}, D_{\bullet}/ D_{\bullet} \cap
V_{\d_1})\big). \end{align*}
Note that $(F'(\j_1,\j_2))_{\j_1,\j_2}$ is a smooth stratification of
$F'$, and that $\kappa_{\j_1,\j_2}$ is a vector bundle. The map $\kappa
q$ decomposes as follows: on $F'(\j_1,\j_2)$, it is equal to the
composition
$$ F'(\j_1,\j_2) \stackrel{\kappa_{\j_1,\j_2}}{\longrightarrow}
\widetilde{\mathcal{F}}_{\j_1}^{im}
\times \widetilde{\mathcal{F}}_{\j_2}^{im}
\stackrel{\pi_1\times \pi_1}{\longrightarrow}
E_{\d_1} \times E_{\d_2}.$$
Similarly,  we define
a smooth stratification $(F''(\i_1, \i_2))_{\i_1,\i_2}$ of $F''$ for
$\i_1, \i_2$ running in the set of complementary subsequences in $\i$,
together with the vector bundles
$\kappa_{\i_1, \i_2}: F''(\i_1,\i_2) \to
\widetilde{\mathcal{F}}_{\i_1}\times \widetilde{\mathcal{F}}_{\i_2}$.
Moreover, one has
$$r^{-1}(F'(\j_1,\j_2))=\bigsqcup_{(\i_1,\i_2)^{im}=(\j_1,\j_2)}
\hspace{-.1in}F''(\i_1,\i_2).$$ Considering the diagram
$$\xymatrix{
 \widetilde{\mathcal{F}}_{\i} \ar[d]^-{\pi_1'}&
\underset{(\i_1,\i_2)^{im}=(\j_1,\j_2)}{\bigsqcup}
\hspace{-.1in}F''(\i_1,\i_2)  \ar[l]_-{\iota''}
\ar[d]^-{r} & &\\
\widetilde{\mathcal{F}}^{im}_{\i}&
F'(\j_1,\j_2) \ar[l]_-{\iota'} \ar[r]^-{\kappa_{\j_1,\j_2}} &
\widetilde{\mathcal{F}}_{\j_1}^{im}
\times \widetilde{\mathcal{F}}_{\j_2}^{im} \ar[r]^-{\pi_1 \times \pi_1}
& E_{\d_1} \times E_{\d_2}}
$$
and reasoning as in \cite[Lemma 9.2.4]{L1}, we see that
$\kappa_{\j_1,\j_2!}(\iota')^* P_{|F'(\j_1,\j_2)}$ belongs to
$\mathcal{Q}_{\j_1} \boxtimes \mathcal{Q}_{\j_2}$ (in particular, it is
semisimple).
Finally, applying \cite[\S 8.1.6]{L1} (see also \cite[\S 4.6, Lemma
4.7]{L2}), we obtain for all $n$
$$H^n((\kappa q)_{!} (\iota')^*P) \simeq \bigoplus_{\j_1,\j_2} H^n
\big(\pi_{1!}\times \pi_{1!}
\kappa_{\j_1,\j_2!} (\iota')^* P_{|F'(\j_1,\j_2)}\big).$$
 Now, $\pi_{1!} \big(\bigoplus_{\j_1,\j_2} \kappa_{\j_1,\j_2!}
(\iota')^*P_{|F'(\j_1,\j_2)}\big)$ and $(\kappa q)_! \iota^*P$ are two
semisimple complexes with isomorphic perverse cohomology, so they are
isomorphic  and the claim follows.~\qed

\vspace{.2in}

\paragraph{\textbf{4.5. The algebra $\mathcal{K}_v$.}} We define an
$\A$-module
$\mathcal{K}(\mathcal{Q}_{\d})$ as follows :
$\mathcal{K}(\mathcal{Q}_{\d})$ is generated by symbols $[P]$ for $P \in
\mathcal{Q}_{\d}$ subject to the relations $[P_1 \oplus P_2]=[P_1] +
[P_2]$ and $[P[n]]=v^n[P]$.
We also set $\mathcal{K}(\mathcal{Q})=\bigoplus_{\d} \mathcal{K}
(\mathcal{Q}_{\d})$, where by convention
$\mathcal{K}(\mathcal{Q}_{0})=\A$. The Verdier duality induces a
semilinear involution $[P] \mapsto
\overline{[P]}=[D(P)]$.
The $\A$-module $\mathcal{K}(\mathcal{Q}_{\d_1} \boxtimes
\mathcal{Q}_{\d_2})$ is defined in a similar fashion, and there is a
canonical isomorphism $\mathcal{K}(\mathcal{Q}_{\d_1} \boxtimes
\mathcal{Q}_{\d_2}) \simeq \mathcal{K}(\mathcal{Q}_{\d_1}) \otimes_{\A}
\mathcal{K}
(\mathcal{Q}_{\d_2})$. Note that the functors $\mathrm{Ind}$ and
$\mathrm{Res}$ commute with direct sums and shifts, and hence induce
$\A$-linear maps $$\mathrm{ind}_{\d_1,\d_2}^{\d_1+\d_2}:
\mathcal{K}(\mathcal{Q}_{\d_1}) \otimes_{\A}
\mathcal{K}(\mathcal{Q}_{\d_2}) \to \mathcal{K}
(\mathcal{Q}_{\d_1+\d_2}),$$
$$\mathrm{res}_{\d_1,\d_2}^{\d_1+\d_2}: \mathcal{K}
(\mathcal{Q}_{\d_1+\d_2})\to \mathcal{K}(\mathcal{Q}_{\d_1})
\otimes_{\A} \mathcal{K}(\mathcal{Q}_{\d_2}).$$

Summing up over $\d_1$ and $\d_2$ yields the maps
$m=\bigoplus_{\d_1,\d_2} \mathrm{ind}_{\d_1,\d_2}^{\d_1+\d_2}:
\mathcal{K}(\mathcal{Q}) \otimes_{\A} \mathcal{K}(\mathcal{Q}) \to
\mathcal{K} (\mathcal{Q})$ and $\Delta'=\bigoplus_{\d_1,\d_2}
\mathrm{res}_{\d_1+\d_2}^{\d_1,\d_2}: \mathcal{K}(\mathcal{Q}) \to
\mathcal{K}(\mathcal{Q}) \otimes_{\A} \mathcal{K}(\mathcal{Q})$.
Finally, we set $\Delta= \barv \circ \Delta' \circ \barv$.

\begin{prop} The space $(\mathcal{K}(\mathcal{Q}),m)$ is an associative
 algebra. Equip the product $\mathcal{K}(\mathcal{Q}) \otimes
\mathcal{K}(\mathcal{Q})$ with a twisted algebra structure by setting
$(x \otimes y)
(x' \otimes y')=v^{(\mathbf{dim}\;x',\mathbf{dim}\;y)}(xx' \otimes
yy')$. Then $\Delta$ is a morphism of algebras.\end{prop}

\noindent
\textit{Proof.} This can be proved exactly in the same way as in
\cite[Chapter 13]{L1}. \qed

\vspace{.1in}

\begin{prop} The set $\{[P]\; P \in \mathcal{P}_{\d}\}$ is an $\A$-basis
of $\mathcal{K}(\mathcal{Q}_{\d})$.\end{prop}

\noindent
\textit{Proof.} By definition, $\{[P]\;|P \in \mathcal{P}_{\d}\}$ is a
generating set of $\mathcal{K}(\mathcal{Q}_{\d})$ over $\A$.
We will show that these elements are linearly independent. We use the
notations and results in the proof of Proposition~4.1.
 Let us call $x \in E_{\d}$
\textit{regular} if for each
$i \in I^{im}$ at least one of the nilpotent elements in $t(x)_i \in
\mathrm{Hom}(V_i,V_i)^{c_i}$ is regular. We denote by $E_{\d}^{reg}$ the
dense open subset of regular elements in $E_{\d}$, and extend this
notation to $\widetilde{\mathcal{F}}_{\i}$ and
$\widetilde{\mathcal{F}}_{\i}^{im}$.
Recall that any
of the simple perverse sheaves in $\mathcal{T}_\i$ is of the form
$IC(X,\mathfrak{L})$ with $X = (u\times Id)^{-1}(Y)$ for a smooth
irreducible  subvariety $Y \subset \mathcal{G}_{\i}^{im}$ and
$\mathfrak{L}=
(u \times Id)^* \mathfrak{K}$ for an irreducible local system
$\mathfrak{K}$ on $Y$.
Put $X^{reg}=X \cap \widetilde{\mathcal{F}}_{\i}^{im,reg}$.

There is a unique complete flag in a vector
space compatible with a given regular nilpotent element. Hence
$\pi_1 : \widetilde{\mathcal{F}}_{\i}^{im,reg} \to \pi_1(
\widetilde{\mathcal{F}}_{\i}^{im,reg})$
is an isomorphism.
Consequently, $\pi_{1!}(IC(X,\mathfrak{L}))=
\pi_{1!}(IC(X^{reg},\mathfrak{L}))=IC(\pi_1(X^{reg},\mathfrak{L})) +
P'$, where $P'$ is a direct sum of simple perverse sheaves supported on
$E_{\d} \backslash E_{\d}^{reg}$.

Finally, assume that $\sum_i \alpha_i [P_i]=0$ is an $\A$-linear
relation between elements $P_i=\pi_{1!}(Q_i) \in \mathcal{P}_{\d}$.
Restricting to the open set $E_{\d}^{reg}$ we deduce a similar linear
relation between the perverse sheaves $Q_i$. But these are simple
perverse sheaves, hence $\A$-linearly independent. Therefore
$\alpha_i=0$ for all $i$ as desired. \qed

\vspace{.1in}

 Define a bilinear form on $\mathcal{K}(\mathcal{Q})$ by the formula
$$([P],[Q])=\sum_{j \in \Z} D_j(D(P),D(Q)) v^j$$
for all $P, Q \in \bigsqcup_{\d} \mathcal{P}_{\d}$. Then $(\,,\,)$ is a
Hopf  pairing, i.e we have
\begin{equation}\label{E:18}
(xy,z)=(x \otimes y, \Delta(z)), \qquad \forall\; x,y,z \in \mathcal{K}
(\mathcal{Q}).
\end{equation}

\vspace{.1in}

Assume $i \in I$ and let $n \geq 1$ with
$n=1$ if $i$ is imaginary. It is easy to see that
the space $E_{n\epsilon_i}$ is a point, and that the
constant sheaf $(\overline{\mathbb{Q}}_l)_{E_{n\epsilon_i}}$
belongs to $\mathcal{P}_{n\epsilon_i}$. The following theorem is
our main result and will be proved in Section~5.

\begin{theo}The assignement
\begin{align*}
F_i^{(n)} &\mapsto [(\overline{\mathbb{Q}}_l)_{E_{n\epsilon_i}}] \qquad
\text{for} \ \ i \in I^{re}, n \geq 1,\\
F_i &\mapsto [(\overline{\mathbb{Q}}_l)_{E_{\epsilon_i}}] \qquad
\text{for} \ \ i \in I^{im}
\end{align*}
extends to an isomorphism of $\A$-algebras $\Theta:
\U^-_{\A}(\g_{A_Q,\mathbf{m}})
\stackrel{\sim}{\to} \mathcal{K}(\mathcal{Q})$ where
$A_Q$ is the Cartan matrix associated
to $Q$ (see Section~2.1) with $m_i=1$ for all $i \in I$.\end{theo}

\vspace{.1in}

\paragraph{\textbf{Definition.}} The set $\mathbf{B}=
\{\Theta^{-1}([P])\;|P \in
\bigsqcup_{\d}\mathcal{P}_{\d}\}$ is
the \textit{canonical basis} of $\U^-_{\A}(\g)$.

\vspace{.1in}

\begin{cor} If $\mathbf{b}, \mathbf{b}' \in \mathbf{B}$, then
$\mathbf{b} \cdot \mathbf{b}'=\sum_{\mathbf{c} \in \mathbf{B}}
a_{\mathbf{c}} \mathbf{c}$ with $a_{\mathbf{c}} \in
\mathbb{N}[v,v^{-1}]$.\end{cor}

\vspace{.2in}

\section{Geometric realization of $\U^-_v(\g)$}

\vspace{.1in}

Note that each imaginary generator $F_{j}$ appears in each of
the defining relations of $\U^-_{\A}(\g)$ with multiplicity at most one.
 In particular, all the relations that we need to check in order to show
that the map $\Theta: \U^-_{\A}(\g_{A_Q,\mathbf{m}})
\to \mathcal{K}(\mathcal{Q})$ is well-defined occur in spaces $E_{\d}$
with $\d_j \in \{0,1\}$ whenever $j \in I^{im}$. But in these cases, all
maps associated to edge loops at an imaginary vertex are necessarily
zero  (by the nilpotency condition) and our construction coincides with
Lusztig's original construction for the quiver obtained from $Q$ by
removing all edge loops. Hence, the relations (\ref{E:04}) follow from
\cite[Theorem 13.2.11]{L1}.

\vspace{.1in}

\paragraph{} It now remains to check that the map $\Theta$ is an
isomorphism. On generators one checks that $\Theta^*(\,,\,)$ and
$\langle\,,\,\rangle$ coincide (up to renormalization).
Using Proposition~4.2 and the fact that both forms are Hopf pairings
(see (\ref{E:18})), we deduce that
$\Theta^*(\,,\,)$ coincides (up to renormalization) with $\langle\,,\,
\rangle$ on $\U^-_{\A}(\g)$. In particular, $\mathrm{Ker}\;\Theta$
belongs to the radical of
 $\langle\,,\, \rangle$ which is trivial by Proposition~1.1. Thus
$\Theta$
is injective.

\vspace{.1in}

In the rest of this section, we prove that $\Theta$ is surjective. We
argue by induction on the dimension vector $\d$. Assume that
$[Q] \in \mathrm{Im}\;\Theta$ for all $Q \in \mathcal{P}_{\d'}$ with
$d'_j \leq d_j$ for all $j$ and $d'_k < d_k$ for at least one $k \in I$.

Let $k \in I$ be a sink (i.e., the only arrows leaving the vertex $k$
are edge loops). If $x \in E_{\d}$, we put
$$n_k(x)=\mathrm{codim}_{V_k}\left(\fq[x_{\sigma_1},\ldots,x_{\sigma_{c_k}}]
\cdot \sum_{j \neq k; j \stackrel{\sigma}{\to} k}
x_{\sigma}(V_j)\right),$$ where $\sigma_1, \ldots, \sigma_{c_k}$ are the
edge loops at $k$, and if $P$ is any complex, we set
$$n_k(P)=\mathrm{inf}_{x \in supp(P)}\; n_k(x).$$

\begin{prop} If $P \in \mathcal{P}_{\d}$ is such that $n_k(P) >0$, then
$[P] \in \mathrm{Im}\;\Theta$.\end{prop}

\noindent
\textit{Proof.} We will prove this by descending induction on $n_k(P)$.
The statement is empty for $n_k(P) \gg 0$. Let us assume that it holds
for all $Q$ with $n_k(Q) >n$ and let us fix some $P$ with $n_k(P)=n$.
Choose
$\i$ and $R \in \mathcal{T}_{\i}$ such that $P=\pi_{1!} R$. We also set
$$\O=\{x \in E_{\d}\;| n_k(x)=n\},$$
$$\O_{\i}=\{(x, D_{\bullet}) \in \widetilde{\mathcal{F}}_{\i}^{im}\;|x
\in \O\}.$$
We will first describe $\mathrm{Res}^{\d}_{\d-n\epsilon_k,
n\epsilon_k}(P)$. For this, fix and embedding $V_{\d-n\epsilon_k}
\subset V_{\d}$ and consider the diagram

$$\xymatrix{
\O_{\i} \ar[d]_{\pi_1} & \O'_{\i} \ar[l]_-{\iota_{\i}}
\ar[d]_{\pi_1} \ar[dr]^p &\\
\O \ar[d]_j & \O' \ar[r]^-{\kappa} \ar[d]_{j'} \ar[l]_-{\iota'}&
\Upsilon_{\d-n\epsilon_k}
\times E_{n\epsilon_k} \ar[d]_{\overline{j}}\\
E_{\d} & F \ar[l]_-{\iota} \ar[r]^-{\kappa}& E_{\d-n\epsilon_k}
\times E_{n\epsilon_k}}
$$
where the notations are as follows : the bottom row is as in
Section~4.3.;
 $\O'=\O \cap F$ and
$\O'_{\i}=\{(x,D_{\bullet}) \in \O_{\i}\;| x \in F\}$;
$$\Upsilon_{\d-n\epsilon_k}=\left\{x \in E_{\d-n\epsilon_k}\;|\;
\fq[x_{\sigma_1},
\ldots x_{\sigma_{c_k}}]\cdot \sum_{j \neq k, j \stackrel{\sigma}{\to}
k} x_{\sigma}(V_j)=(V_{\d-n\epsilon_k})_k\right\}$$
and all maps are obvious ones (with $p=\kappa\pi_1$). Note that all
squares are cartesian. Let $l: \O_{\i} \to
\widetilde{\mathcal{F}}_{\i}^{im}$ be the inclusion. We have $j^*
\pi_{1!}(R)=\pi_{1!}l^*(R)$. Moreover, $l^*(R)$ is perverse since
$supp(R) \subset \overline{\O_{\i}}$ and $supp(R) \cap \O_{\i} \neq
\emptyset$. On the other hand,
\begin{equation*}
\begin{split}
\overline{j}^* \kappa_!\iota^*\pi_{1!}(R)&=
\kappa_!(j')^*\iota^*\pi_{1!}(R)\\ &= \kappa_! (\iota')^*j^*\pi_!(R)\\
&= \kappa_! (\iota')^* \pi_{1!} l^*(R)\\
&= \kappa_!\pi_{1!} \iota_{\i}^*l^*(R)\\
&=p_! \iota_{\i}^*l^*(R).
\end{split}
\end{equation*}
Let us denote for simplicity $\iota_{\i}^*l^*(R)$ by $R'$.
Let $H$ be the parabolic subgroup of $G_{\d}$ stabilizing
$V_{\d-n\epsilon_k}$. We have $\O_{\i}=G_{\d} \times_H \O'_{\i}$ so that
$(\iota')^*: \mathcal{M}_{G_{\d}} (\O_{\i}) \to
\mathcal{M}_H(\O'_{\i})[\mathrm{dim}\;G_{\d}/H]$ is an  equivalence of
categories. In particular, $R'[-\mathrm{dim}\;G/H]$ is perverse.

\vspace{.1in}

Set $\O^{reg}=
\O \cap E_{\d}^{reg}$ and $\O_{\i}^{reg}=(\pi_1)^{-1}(\O^{reg})$.
Note that from the proof of Proposition~4.3 it follows that
$supp(R)=X=(u \times Id)^{-1}(Y)$ for some irreducible $Y \subset
\mathcal{G}_{\i}^{im}$. Then $X^{reg}$ is open in $X$ and thus $X^{reg}
\cap \mathcal{O}_{\i}$ is open in $X$. Thus $R=IC(X^{reg} \cap
\mathcal{O}_{\i})$ hence $supp(l^*R) \subset supp(R) \subset
\overline{\O_{\i}^{reg}}$. Now, if $(x,D_{\bullet}) \in \O_{\i}^{reg}$
then $D_{\bullet}\cap V_k$ is completely determined by $x$. On the other
hand,  the subspace
$\fq[x_{\sigma_1},
\ldots x_{\sigma_{c_k}}]\cdot \sum_{j \neq k, j \stackrel{\sigma}{\to}
k} x_{\sigma}(V_j)$ is stable under $x_{\sigma_i}$ for $i=1, \ldots ,
c_k$, so it has to belong to $D_{\bullet}\cap V_k$. Define
a closed subset of $\O_{\i}(\i)\subset\O_{\i}$ by the condition
$$\fq[x_{\sigma_1},
\ldots x_{\sigma_{c_k}}]\cdot \sum_{j \neq k, j \stackrel{\sigma}{\to}
k} x_{\sigma}(V_j)\in D_{\bullet} \cap V_k.\qquad \qquad\qquad (*)$$
Then we have $supp(l^*R) \subset \O_{\i}(\i)$. We define
$\O'_{\i}(\i)$ similarly
and we have $supp(R') \subset \O'_{\i}(\i)$. Construct a sequence
$\j$ by deleting the last $d_k-n$ entries equal to $k$ in $\i$ (so that
if $(x,D_{\bullet}) \in \O'_{\i}(\i)$, then $D_{\bullet} \cap
V_{\d-n\epsilon_k}$ is of type $\j$). Now consider the diagram

$$\xymatrix{
\widetilde{\O}_{\i}(\i) \ar[r]^-{\kappa''} \ar[d]_{\pi_1'} &
\Upsilon_{\j} \times \widetilde{\mathcal{F}}_{n\epsilon_k}
\ar[d]_{\pi_1'} \ar[r]^-{h'}&  \widetilde{\mathcal{F}}_{\j} \times
\widetilde{\mathcal{F}}_{n\epsilon_k} \ar[d]_{\pi_1'}\\ \O'_{\i}(\i)
\ar[r]^-{\kappa'} \ar[d]_{\pi_1} & \Upsilon^{im}_{\j} \times
\widetilde{\mathcal{F}}^{im}_{n\epsilon_k} \ar[d]_{\pi_1} \ar[r]^-{h}&
\widetilde{\mathcal{F}}^{im}_{\j} \times
\widetilde{\mathcal{F}}^{im}_{n\epsilon_k}  \ar[d]_{\pi_1}\\
\O' \ar[r]^-{\kappa} & \Upsilon_{\d-n\epsilon_k} \times E_{n\epsilon_k}
\ar[r]^{\overline{j}} & E_{\d-n\epsilon_k} \times E_{n\epsilon_k}
}
$$
with $\widetilde{\O}_{\i}(\i)=(\pi_1')^{-1}(\O'_{\i}(\i))$,
$\Upsilon^{im}_{\j}=\pi_1^{-1}(\Upsilon_{\d-n\epsilon_k}) \subset
\widetilde{\mathcal{F}}^{im}_{\j}$, $\Upsilon_{\j}=(\pi_1')^{-1}
(\Upsilon^{im}_{\j})$ and with obviously defined maps ($h$ and $h'$ are
open embeddings).  Note that
$\kappa, \kappa', \kappa''$ are all vector bundles of rank $c_k(d_k-n)n$
by the condition $(*)$ together with the fact that $k$ is a sink. From
this one deduces that the two squares in the above diagram are
cartesian. By construction, $R'$ appears in $\pi_{1!}'(
(\overline{\mathbb{Q}_l})_{\widetilde{\O}_{\i}(\i)})=\pi'_{1!}
(\kappa'')^*((\overline{\mathbb{Q}_l})_{\Upsilon_{\j} \times
\widetilde{\mathcal{F}}_{n\epsilon_k}})=(\kappa')^*
\pi'_{1!}((\overline{\mathbb{Q}_l})_{\Upsilon_{\j} \times
\widetilde{\mathcal{F}}_{n\epsilon_k}})=(\kappa')^*h^* \pi'_{1!}
((\overline{\mathbb{Q}_l})_{\widetilde{\mathcal{F}}_{\j} \times
\widetilde{\mathcal{F}}_{n\epsilon_k}})$.
In particular, there exists $R'' \in
\mathcal{M}(\widetilde{\mathcal{F}}^{im}_{\j} \times
\widetilde{\mathcal{F}}^{im}_{n\epsilon_k})[\mathrm{dim}\;G_{\d}/H +
c_k(d_k-n)n]$ such  that $R'=(\kappa')^*
h^*R''$.
Hence $\pi_{1!}R'=\kappa^* \pi_{1!}h^*(R'')$
and altogether we obtain
$$\overline{j}^*\widetilde{\mathrm{Res}}^{\d}_{\d-n\epsilon_k,
n\epsilon_k}(P)= p_!(R')
=\kappa_!\kappa^* \pi_{1!}h^*(R'')=\pi_{1!}h^*(R'')[-2c_k(d_k-n)n],$$
which is an object in
$\mathcal{M}(\Upsilon_{\d-n\epsilon_k}\times
E_{n\epsilon_k})[\mathrm{dim}\; G_{\d}/H-c_k(d_k-n)n]$ by
Proposition~4.1.
Further, note that $\pi_{1!}h^*(R'')=\overline{j}^*\pi_{1!}(R'')$ and
$\pi_{1!}(R'') \in \mathcal{Q}_{\d-n\epsilon_k}
\boxtimes \mathcal{Q}_{n\epsilon_k}$ by construction. Hence
$[\pi_{1!}(R'')]  \in \mathrm{Im}\;\Theta \otimes \mathrm{Im}\;\Theta$
by our first induction  hypothesis. On the other hand, all squares in
the following diagram are cartesian~: $$\xymatrix{
E_{\d} &E'' \ar[l]_-{p_3} & E' \ar[l]_-{p_2} \ar[r]^-{p_1}&
E_{\d-n\epsilon_k} \times E_{n\epsilon_k}\\
\O' \ar[u]^-{j\iota'} & \O' \ar[u] \ar[l]_{Id} & \O' \times
G_{\d-n\epsilon_k} \times GL(n) \ar[u] \ar[l] \ar[r] &
\Upsilon_{\d-n\epsilon_k} \times E_{n\epsilon_k} \ar[u]^{\overline{j}}\\
& & \O' \ar[ul]_-{Id} \ar[u] \ar[ur]^-{\kappa} &}$$

It then follows from the definitions that
\begin{equation*}
\begin{split}
(\iota')^*j^*
\widetilde{\mathrm{Ind}}_{\d-n\epsilon_k,n\epsilon_k}^{\d}(\pi_{1!}
R'')&=\kappa^* \overline{j}^* \pi_{1!}(R'')\\
&=\kappa^*\pi_{1!}h^*(R'')\\
&=\pi_{1!}(R')\\
&=(\iota')^*j^*(P).
\end{split}
\end{equation*}
Since $(\iota')^*:\mathcal{M}_{G_{\d}}
(\O_{\i}) \to \mathcal{M}_H(\O'_{\i})[\mathrm{dim}\;G_{\d}/H]$ is fully
faithful, we deduce \begin{equation}\label{E:666}
j^*(P)=j^*\widetilde{\mathrm{Ind}}_{\d-n\epsilon_k,n\epsilon_k}^{\d}(\pi_{1!}
R'').
\end{equation}
But it is obvious that
$supp(\widetilde{\mathrm{Ind}}_{\d-n\epsilon_k,n\epsilon_k}^{\d}(\pi_{1!}
R'')) \in \overline{\O}$. Thus from (\ref{E:666}) we deduce that,  in
$\mathcal{Q}_{\d}$, $[P]=
[\widetilde{\mathrm{Ind}}_{\d-n\epsilon_k,n\epsilon_k}^{\d}(\pi_{1!}
R'')]+[P']$, where $P'$ has support in $\overline{\O}\backslash \O$.
Such  complex necessarily satisfies $n_k(P') >n$. Therefore, by our
second induction hypothesis $[P'] \in \mathrm{Im}\;\Theta$ and thus
$[P] \in \mathrm{Im}\;\Theta$. This completes the proof of
Proposition~5.1 \qed

\vspace{.1in}

\paragraph{}The rest of the proof of Theorem~4.1 goes as in \cite{L1}.
Namely, let $P \in \mathcal{P}_{\i}$ with $\i=(i_1, \ldots, i_l,k)$.
Using the Fourier-Deligne transform, we may assume that $k$ is a sink in
our quiver. But  then it is easy to see that $n_k(P)>0$ and from
Proposition~5.1 it follows that $[P] \in \mathrm{Im}\;\Theta$ as
desired.\qed

\vspace{.1in}

\begin{cor}\label{C:1} \hfill

{\rm (i)} Fix $k \in I$ and $n \geq 1$. Set $F_k^{(n)}=F_k^n$ if
$k \in I^{im}$.
Then there exist subsets
$\mathbf{B}_{\geq n,k}, \mathbf{B}'_{\geq n,k} \subset \mathbf{B}$ such
that $$
F_k^{(n)} \mathbf{U}_{\mathbb{A}}^-=\bigoplus_{\mathbf{c} \in
\mathbf{B}_{\geq n,k}}
\mathbb{A} \mathbf{c}, \qquad  \mathbf{U}_{\mathbb{A}}^-F_k^{(n)}=
\bigoplus_{\mathbf{c} \in \mathbf{B}'_{\geq n,i}}\mathbb{A}
\mathbf{c}.$$

{\rm (ii)} For every $b \in \mathbf{B}$,
there exists $k$ such that $b \in \mathbf{B}_{\geq 1, k}$.\\

{\rm (iii)} If $b \in \mathbf{B}_{\geq n, k}\backslash \mathbf{B}_{\geq
n+1, k}$, then there exists
$b' \in \mathbf{B}$ such that
$$F_k^{(n)}b' \in b \oplus \bigoplus_{\mathbf{c} \in \mathbf{B}_{\geq
n+1,k}} \A \mathbf{c}.$$
\end{cor}

\noindent
\textit{Proof.} We use the same notations as in the proof of
Theorem~4.1.  In (i), we only prove the statement concerning the space
$\mathbf{U}_{\mathbb{A}}^-F_k^{(n)}$. The other part is proved in an
analogous way. Note that $\mathbf{P}_{n\epsilon_k}=\{P_{n,k}\}$ consists
of a single element.
Let $x \in \U_{\mathbb{A}}^-[\d]$.
Write $\Theta(x)=\sum \alpha_P
[P]$, so that $\Theta(x F_k^{(n)})=\sum \alpha_P
[\mathrm{Ind}_{\d,n\epsilon_k}^{\d+n\epsilon_k}(P \boxtimes P_{n,k})]$.
By construction, it is clear that
$n_k(\mathrm{Ind}_{\d,n\epsilon_k}^{\d+n\epsilon_k}(P \boxtimes
P_{n,k})) \geq n$ for any $P$. Thus
$$\Theta(\U_{\A}^- F_i^{(n)}) \subset \bigoplus_{n_k(Q) \geq n} \A
[Q].$$ We now prove the opposite inclusion. Fix $\d \in \N^I$. We argue
by descending induction on $n_k(Q)$. The statement is empty for $n_k(Q)
> d_k$. Fix $R \in \mathcal{P}_{\d}$ with $n_k(P) >0$ and assume that
$[Q] \in \Theta(\U_{\A}^- F_k^{(n_k(Q))})$ for all
$Q \in \mathcal{P}_{\d}$
such that $n_k(Q) > n_k(R)$.
The proof of Theorem~4.1 shows the existence of $T \in \mathcal{P}$ such
that \begin{equation}\label{E:enddas}
[R]=[T] \cdot [P_{n_k(R),k}] + [P']\qquad \text{with}\;n_k(P') > n_k(R).
\end{equation}
By induction
$[P'] \in \Theta (\U_{\A}^- F_i^{(n_k(P'))}) \subset
\Theta (\U_{\A}^- F_i^{(n_k(R))})$ and we are done.

The statement (ii) follows from the Fourier-Deligne transform (see the
end of the proof of Theorem~4.1). The statement (iii) follows from
(\ref{E:enddas}) and the dual version  of (\ref{E:enddas}).
\qed

\vspace{.1in}

\begin{cor}\label{C:2} For any $b \in \mathbf{B}$ we have
$\overline{b}=b$. Assume that $c_i >1$ for all $i \in I^{im}$. Then the
following holds : \begin{equation}\label{E:neqqq1}
\forall \; b \in \mathbf{B} \qquad \langle b,b \rangle \in 1 + v
\Z[[v]], \end{equation}
\begin{equation}\label{E:neqqq2}
\forall\; b \neq b' \in \mathbf{B} \qquad \langle b, b' \rangle \in v
\Z[[v]]. \end{equation}
\end{cor}
\noindent
\textit{Proof.} The first statement can proved
in the same way as in \cite[\S 13]{GL}.
When $c_i > 1$ for all $i \in I^{im}$ the elements of
$\sqcup_{\d}\mathcal{P}_{\d}$ are all simple perverse sheaves by
Proposition~4.1.
The relations (\ref{E:neqqq1}) and (\ref{E:neqqq2})
then follow from \cite[\S 8.1.10]{L1}
together with the fact that the forms $\Theta^*(\,,\,)$ and
$\langle\;, \; \rangle$ differ (on each weight space) by a factor in $1
+ v \Z[[v]]$.\qed

\vspace{.1in}

\begin{prop}\label{C:3} Let $\mathcal{B}$ be the global basis of
$\U_{\A}^-(\g)$ defined in \cite{JKK}. Following \cite{JKK} we denote
the elements of $\mathcal{B}$ by $G(\beta)$.  Then the following
statements hold\,:
\begin{enumerate}
\item[i)] For any $G(\beta) \in \mathcal{B}$ we have
$\overline{G(\beta)}=G(\beta)$. \item[ii)] Assume that $c_i > 1$ for all
$i \in I^{im}$. Then we have \begin{equation*}
\forall \; G(\beta) \in \mathcal{B} \qquad \langle G(\beta),G(\beta)
\rangle \in 1 + v \Z[[v]], \end{equation*}
\begin{equation*}
\forall\; G(\beta) \neq G(\beta') \in \mathcal{B} \qquad \langle
G(\beta), G(\beta') \rangle \in v \Z[[v]]. \end{equation*}
\item[iii)]For every $i \in I$ and $n \geq 1$,
there exists a subset $\mathcal{B}_{\geq n,i} \subset \mathcal{B}$ such
that
$$F_i^{(n)}\mathbf{U}_{\mathbb{A}}^-\bigoplus_{G(\gamma) \in
\mathcal{B}_{\geq n,i}}
\mathbb{A} G(\gamma).$$
\item[iv)] For any $G(\beta) \in  \mathcal{B}_{\geq n, i}\backslash
\mathcal{B}_{\geq n+1, i}$ there exists $G(\beta') \in \mathcal{B}$ such
that
$$F_i^{(n)}G(\beta') \in G(\beta) \oplus
\bigoplus_{G(\gamma) \in \mathcal{B}_{\geq n+1,i}} \A G(\gamma).$$
\end{enumerate}
\end{prop}
\noindent
\textit{Proof.} The first three statements can be found in \cite{JKK},
in  Sections~ 9.3, 7.38 and 10.2, respectively.
We prove the last statement. Let $(L_{\infty}, \mathcal{B}_{\infty})$
denote the crystal basis of $\U^-_{\A}(\g)$ defined in \cite{JKK} and
let $\tilde{e}_i, \tilde{f}_i, e'_i$ stand for the Kashiwara operators.
Let $G(\beta) \in  \mathcal{B}_{\geq n, i}\backslash \mathcal{B}_{\geq
n+1, i}$. Denote by $\beta \in \mathcal{B}_{\infty}$ the corresponding
element in the crystal graph. Then there exists $\beta_0$ in
$\mathcal{B}_{\infty}$ such that
$\tilde{e}_i \overline{{\beta}_0}=0$ and $\tilde{f}_i^n
\overline{\beta_0}= \beta$. Let $G(\beta_0) \in \mathcal{B}\backslash
\mathcal{B}_{\geq 1, i}$ be the global basis element associated to
$\beta_0$. By $iii)$, there exists $G(\beta_i) \in \mathcal{B}$ and $a_i
\in \A$ such that $$F_i^{(n)}G(\beta_0) \equiv \sum_i a_i
G(\beta_i)\;(\mathrm{mod}\;F_i^{(n+1)}\U_{\A}^-(\g)).$$ Moreover we have
$\overline{a_i}=a_i$ by invariance of $G(\beta_0),G(\beta_i)$ under the
bar involution. Let $P_i$ denote the projection of
$\U_{\A}^-(\g)=\mathrm{Ker}\; e'_i \oplus F_i \U_{\A}^-(\g)$ onto
$\mathrm{Ker}\;e'_i$. We have
\begin{alignat*}{2}
F_i^{(n)}G(\beta_0) & \equiv
F_i^{(n)}P_iG(\beta_0)\;&(\mathrm{mod}\;F_i^{(n+1)}\U_{\A}^-(\g))\\ &
\equiv \tilde{f}_i^n P_i
G(\beta_0)\;&(\mathrm{mod}\;F_i^{(n+1)}\U_{\A}^-(\g))\\ & \equiv
\tilde{f}_i^n G(\beta_0)\;&(\mathrm{mod}\;F_i^{(n+1)}\U_{\A}^-(\g))\\ &
\equiv G(\beta)\;&(\mathrm{mod}\;F_i^{(n+1)}\U_{\A}^-(\g)\oplus
vL_{\infty}) \end{alignat*}
Thus we deduce that $a_i\equiv 1 \;(\mathrm{mod}\;v\Z[v])$ if
$\beta_i=\beta$ and $a_i \in v \Z[v]$ otherwise. But from
$a_i=\overline{a_i}$ it now follows that $a_i=1$ for $\beta_i=\beta$ and
$a_i=0$ otherwise.\qed

\vspace{.1in}

>From Corollary~\ref{C:1}, Corollary~\ref{C:2} and Proposition~\ref{C:3}
> we deduce, using the method in
\cite{GL}, the following result.
\begin{theo}\label{T:1}If $c_i > 1$ for all $i \in I^{im}$ then
$\mathbf{B}$ coincides with the global basis $\mathcal{B}$.\end{theo}

\vspace{.1in}

We conjecture that Theorem~\ref{T:1} holds unconditionally.

\vspace{.2in}

\small{

}




\end{document}